\documentclass[11pt]{amsart}
\usepackage[utf8x]{inputenc}
\usepackage{graphicx, xcolor}
\usepackage{amsmath,amssymb,amsthm,comment,mathtools}
\usepackage{wrapfig}
\usepackage[colorlinks,allcolors=blue]{hyperref} 
\usepackage{cleveref}
\usepackage{pdfpages}
\usepackage{amsmath}
\usepackage{subfigure,cite}
\usepackage{tikz}
\usepackage[a4paper,margin=2.6cm]{geometry}

\theoremstyle{plain}
\newtheorem{thm}{Theorem}[section]

\newtheorem{prop}[thm]{Proposition}
\newtheorem{lem}[thm]{Lemma}

\newtheorem{introtheorem}{Theorem}

\newenvironment{theoremproof}[1]
  {\begin{proof}[Proof of Theorem~\ref{#1}]}
  {\end{proof}}

\theoremstyle{definition}
\newtheorem{ex}[thm]{Example}
\newtheorem{defi}[thm]{Definition}
\newtheorem{rmk}[thm]{Remark}

\newcommand{\Hess}{\mathcal{H}\mathrm{ess}}
\newcommand{\mx}{\mathsf{X}}

\newcommand{\hess}{{\mathbf h}}

\newcommand{\Fl}{\mathcal{F}\ell ags}

\newcommand{\C}{\mathbb{C}}           
\newcommand{\Z}{\mathbb{ Z}}           
 
\newcommand{\Q}{\mathbb{ Q}}           

\newcommand{\ForcedFirstColumn}[1]{%
  \begingroup
  \pgfmathtruncatemacro{\pdq}{#1}
  \pgfmathtruncatemacro{\pdlast}{#1-1}
  \begin{tikzpicture}[
    x=9mm, y=-9mm,
    baseline=(current bounding box.center),
    font=\small
  ]
    \fill[blue!6] (0,0) rectangle (1,\pdlast);

    \draw[gray!60] (0,0) rectangle (1,\pdq);
    \foreach \pdrow in {1,...,\pdq} {
      \draw[gray!60] (0,\pdrow) -- (1,\pdrow);
      \node[anchor=east] at (-0.4,{\pdrow-0.5})
        {$\pdrow$};
    }

    \ifnum\pdlast>0
      \foreach \pdrow in {1,...,\pdlast} {
        \draw[line width=0.8pt]
          (0,{\pdrow-0.5}) -- (1,{\pdrow-0.5});
        \node[anchor=west] at (1.1,{\pdrow-0.5})
          {$x_{\pdrow}$};
      }
    \fi

    \draw[line width=0.8pt,rounded corners=4pt]
      (0.5,\pdq) -- (0.5,{\pdq-0.5})
      -- (1,{\pdq-0.5});

    \draw[blue!70!black,line width=1.3pt,
          rounded corners=4pt,->]
      (-0.3,{\pdq-0.5}) -- (0.5,{\pdq-0.5})
      -- (0.5,-0.3);

    \node at (0.5,-0.65) {$1$};
    \node[anchor=west] at (1.1,{\pdq-0.5}) {elbow};
    \node at (0.5,{\pdq+0.65})
      {Forced weight: $\displaystyle\prod_{r=1}^{\pdlast}x_r$};
  \end{tikzpicture}%
  \endgroup
}

\usepackage{tikz}
\usetikzlibrary{calc}

\tikzset{
  pdpipe/.style={line width=1pt,line cap=round,line join=round},
  pdgrid/.style={draw=black!20,line width=.3pt}
}

\newcommand{\exampleqed}{%
    {   
    \renewcommand{\qedsymbol}{$\lozenge$}%
    \qed%
    }
}

\numberwithin{equation}{section}
\title[Presentations of rank 1 nilpotent Hessenberg varieties]
{Presentations of the cohomology ring of 
Hessenberg varieties for the minimal nilpotent orbit in type A}

\author[G. Gaggero, et al]{Giulia Gaggero}
\address[G. Gaggero] {Department of Mathematics and Statistics
McMaster University, Hamilton, ON L8S 4L8, Canada}
\email{gaggerog@mcmaster.ca}

\author[]{Megumi Harada}
\address[M. Harada] {Department of Mathematics and Statistics
McMaster University, Hamilton, ON L8S 4L8, Canada}
\email{haradam@mcmaster.ca}

\author[]{Hannah Nardone}
\address[H. Nardone]{Department of Mathematics, McGill University,
Montreal, QC
H3A 0B9, Canada} 
\email{hannah.nardone@mail.mcgill.ca}

\author[]{Brett Nasserden}
\address[B. Nasserden] {Department of Mathematics and Statistics
McMaster University, Hamilton, ON L8S 4L8, Canada}
\email{nasserdb@mcmaster.ca }

\author[]{Jenna Rajchgot}
\address[J. Rajchgot]
{Department of Mathematics and Statistics
McMaster University, Hamilton, ON L8S 4L8, Canada}
\email{rajchgoj@mcmaster.ca }

\author[]{Lisa Seccia}
\address[L. Seccia]
{Department of Mathematics and Statistics
McMaster University, Hamilton, ON L8S 4L8, Canada}
\email{seccial@mcmaster.ca }

\author[]{Adam Van Tuyl}
\address[A. Van Tuyl]{Department of Mathematics and Statistics
McMaster University, Hamilton, ON L8S 4L8, Canada}
\email{vantuyla@mcmaster.ca}

\date{\today}

\begin{document}

\begin{abstract} 
In 2016 Abe and Crooks found a presentation for the cohomology
ring of Hessenberg varieties associated with the minimal nilpotent orbit. 
Specifically, they presented the ring as a quotient of the cohomology ring of the flag variety, and described the kernel of the quotient map as the span of certain Schubert classes. This note improves upon their 
result in the Lie type A case, by giving a minimal generating set for the kernel. Our description uses substantially fewer generators than that of Abe-Crooks. Our list consists of a set of 
Schubert polynomials associated to permutations constructed
from the Hessenberg function $\hess$.
\end{abstract}

\maketitle

\section{Introduction}\label{sec: intro} 

Let $n$ be a positive integer and let $\Fl(\C^n)$ denote the flag variety of nested sequences of subspaces $V_\bullet=(0\subset V_1\subset V_2 \subset \cdots \subset V_n = \C^n)$ such that $\dim (V_i)=i$ for all $i$. Such a sequence $V_\bullet$ of subspaces is called a \emph{full flag}. Let $\hess: [n]:=\{1,2,\cdots,n\} \to [n]$ be a weakly increasing function such that $\hess(i)\geq i$ for all $i$, called a \emph{Hessenberg function}. The subvariety 
$$
\Hess(\mx, \hess) := \{V_\bullet \in \Fl(\C^n): \mx (V_i)\subseteq V_{\hess(i)} \, \textup{ for all } \, i \in [n] \}
$$
is the (type $A$) \emph{Hessenberg variety} in $\Fl(\C^n)$ associated to a choice of Hessenberg function $\hess$ and a matrix $\mx \in \mathfrak{gl}_n(\C)$. The study of Hessenberg varieties lies at the intersection of algebraic geometry, representation theory, combinatorics, and theory of quasisymmetric functions, 
among many other research areas (cf. \cite[Part III]{Handbook-Hess} and references therein). In this manuscript we focus exclusively on Lie type A, i.e., the flag variety corresponding to $GL_n(\C)$.

The study of the cohomology rings of Hessenberg varieties has been an active research area for some time. We view this theory as part of a rich history of the study of the cohomology rings of subvarieties of (partial and full) flag varieties (as a non-exhaustive list, see e.g. \cite{GasharovReiner, MbirikaTymoczko, ReinerWooYong, StDizierYong}).  We also note that a recent preprint \cite{AGR2026} studies a situation very similar to ours, namely, the cohomology rings of certain unions of Schubert varieties. 
One of the challenges of this research area within Hessenberg geometry is to find explicit descriptions of these cohomology rings via generators and relations, in a manner similar to the well-known case of $\Fl(\C^n)$ with its famous Borel presentation. When $\mx$ is nilpotent, the case $\hess(i)=i$ for all $i$ is a well-studied special case where the corresponding Hessenberg varieties are called \textit{Springer fibers}. The cohomology rings of Springer fibers were studied several decades ago (e.g. De Concini--Procesi \cite{DeConciniProcesi}, Garsia--Procesi \cite{GarsiaProcesi}) and have an explicit description via the Tanisaki ideal 
(e.g. \cite{GarsiaProcesi, Tanisaki}). However, outside of the 
Springer case, less is known, although some cases are understood (e.g. \cite{AHHM2019, AHM2019}).

Abe and Crooks study Hessenberg varieties associated with the minimal nilpotent orbit in \cite{AbeCrooks}. Specifically, they give a description of the cohomology ring of the associated Hessenberg variety $\Hess(\mx, \hess)$ as a quotient of the cohomology ring of the flag variety, where the kernel of the quotient map is additively generated by an explicit list of Schubert classes \cite[Theorem 32]{AbeCrooks}. In the present manuscript, we improve upon their result by finding a substantially smaller set of Schubert classes which minimally generate this kernel (Theorem~\ref{thm:intro-cohomology}). 

It should be acknowledged that there has been other work in Hessenberg geometry which has built on, or generalized, results in Abe and Crooks' work which is not discussed in this manuscript. Examples include Goldin and Precup's flat degenerations from minimal semisimple to minimal nilpotent Hessenberg schemes \cite{GoldinPrecup-flat}, Escobar, Precup, and Shareshian's results on singular loci, irreducibility, reducedness, and Poincaré polynomials of  codimension one Hessenberg varieties in type A \cite{EscobarPrecupShareshian}, and Goldin and Tymoczko's work on the GKM properties enjoyed by certain Hessenberg varieties \cite{GoldinTymoczko}.

We work with $\Q$ coefficients for our cohomology rings. We also assume throughout the paper that $n$ is a positive integer, $n \geq 2$. Here and below, we fix once and for all our choice $\mx = E_{1,n}$, where $E_{1,n}$ is the $n\times n$ matrix with a $1$ in position $(1,n)$ and a $0$ elsewhere. In type A, this is a minimal nilpotent element, and the results of \cite{AbeCrooks} apply. With this established, we recall the following.

\begin{defi}\label{def:forbidden}
Let $w \in S_n$ and let $\hess: [n] \to [n]$ be a Hessenberg function. We say that $w$ is $\hess$-forbidden if $w^{-1}(1)>\hess(w^{-1}(n))$.  We let $\mathcal{F}_\hess$ denote the set of $\hess$-forbidden permutations. 
\end{defi}

The above terminology is motivated by a geometric result in the work of Abe and Crooks: they show in \cite[Proposition 3]{AbeCrooks} that $w$ is $\hess$-forbidden precisely when its associated coordinate flag does not lie in $\Hess(\mx = E_{1,n},\hess)$. They then show that the cohomology ring $H^*(\Hess(\mx,\hess);\Q)$ is isomorphic to $H^*(\Fl(\C^n);\Q)/J_\hess$ where $J_\hess := \langle \sigma_w: w \in \mathcal{F}_\hess\rangle$ and $\sigma_w$ denotes the Schubert class corresponding to $w$. 

The main results of this note give a substantially smaller generating set for the denominator of the above description of $H^*(\Hess(\mx,\hess);\Q)$. In fact, we give two separate results, as we now explain. Recall that the cohomology of the flag variety $H^*(\Fl(\C^n);\Q)$ has a Borel description as $\Q[x_1,\cdots,x_n]/\langle \mathbf{e}_1,\ldots, \mathbf{e}_n \rangle $ as in \cite{Borel1953}, where $\mathbf{e}_i$ is the $i^{\textnormal{th}}$ elementary symmetric polynomial in $x_1,\ldots, x_n$. In fact, this description is also valid with $\Z$ coefficients. 
For $w\in S_n$, let $\mathfrak{S}_w$ denote the associated Schubert polynomial~\cite{LascouxSchutzenberger1982} in the variables $x_1,\ldots, x_n$. It is well-known that the Schubert cohomology class $\sigma_w\in H^{2\ell(w)}(\Fl(\C^n);\Q)$ of the opposite Schubert variety associated to $w$ can be represented by $[\mathfrak{S}_w]$ in the above Borel presentation; we use this (opposite) convention for Schubert classes throughout (cf. \cite[Part I, Chapter 1]{Handbook-Hess} and \cite{AbeCrooks}). In analogy to the ideal $J_\hess$ above, we now define $\tilde{J}_\hess := \langle \mathfrak{S}_w: w \in \mathcal{F}_\hess\rangle \subset \Z[x_1,\cdots,x_n]$. Our two main results, Theorem~\ref{thm:intro-polynomial} and Theorem~\ref{thm:intro-cohomology}, give efficient sets of generators of $\tilde{J}_\hess$ and $J_\hess$ respectively.  

We take a moment to discuss the motivation for presenting our results in the form given below. It would have been possible to present Theorem~\ref{thm:intro-cohomology} as our one main result. We chose instead to present two separate theorems because it makes it clearer that the simplification of the set of generators, which was the original motivation for our investigations, can be understood in two stages. The simplification that is recorded in Theorem~\ref{thm:intro-polynomial}(2) occurs at the level of the polynomial ring, arising from a relationship among Schubert polynomials. It is when we quotient by the ideal generated by positive degree symmetric polynomials (to obtain the cohomology of the flag variety) that we obtain the next simplification, recorded in Theorem~\ref{thm:intro-cohomology}, which allows us to also delete some of the Schubert classes $\sigma_{c_j^{(\hess)}}$ from the list of generators. 

We now define the set of permutations that give our efficient generating sets. 
For the definition below, we represent permutations by their one-line notation, so for $w \in S_n$ we write $w = (w(1),\cdots,w(n))$ and refer to $w(j)$ as the entry in position $j$. 

\begin{defi}\label{def:gens}
Let $\hess:[n]\to[n]$ be a Hessenberg function. For each $j\in[n]$
with $\hess(j)<n$, let $c_j^{(\hess)}\in S_n$ be the permutation with $n$ in
position $j$, $1$ in position $\hess(j)+1$, and all remaining entries arranged 
in increasing order. In one-line notation, we have
\[
c_j^{(\hess)}=(2,\ldots,j,n,j+1,\ldots,\hess(j),1,\hess(j)+1,\ldots,n-1),
\]
where empty intervals are assumed to be omitted.
\end{defi}

 Our first result gives a monic Gr\"obner basis for $\tilde{J}_\hess:= \langle \mathfrak{S}_w: w \in \mathcal{F}_\hess\rangle$ in $\Z[x_1,\cdots,x_n]$. 

\begin{introtheorem}\label{thm:intro-polynomial}
Let $n\geq 2$, $\hess:[n]\to[n]$ be a Hessenberg function.   
Then: 
\begin{enumerate} 
\item for each $j \in [n], \hess(j)<n$, we have 
\[
\mathfrak{S}_{c_j^{(\hess)}}
=
\left(\prod_{r=1}^{\hess(j)}x_r\right)
H_{n-j-1}(x_1,\ldots,x_j)
\]
where $H_d(x_1,\ldots,x_j)$ denotes the complete homogeneous
symmetric polynomial of degree $d\geq 0, d\in \Z$ in the variables $x_1,\ldots,x_j$,
with $H_0 :=1$, 
\item as ideals in $\mathbb{Z}[x_1,\ldots,x_n]$ we have 
\[
\tilde{J}_\hess := \bigl\langle \mathfrak{S}_w : w\in\mathcal{F}_\hess \bigr\rangle
=
\bigl\langle \mathfrak{S}_{c_j^{(\hess)}} : j\in[n],\ \hess(j)<n \bigr\rangle
\]
\item for the lexicographic order $x_n>\cdots>x_1$, the polynomials $\{\mathfrak{S}_{c_j^{(\hess)}}: j \in [n], \hess(j)<n\}$  
form a monic Gr\"obner basis for $\tilde{J}_\hess$. 
\end{enumerate}  
\end{introtheorem}

We note that the generating set in Theorem~\ref{thm:intro-polynomial} need not be minimal (see Remark
\ref{rmk: redundant generators}).
The above theorem shows that the set $\{\mathfrak{S}_{c_j^{(\hess)}}: j \in [n], \hess(j) <n\}$ is sufficient to generate $\tilde{J}_\hess$. The analogous result for the ideal $J_\hess$ in $H^*(\Fl(\mathbb{C}^n));\Q)$ requires even fewer generators, indexed by a set $D(\hess)$ defined below. We have the following. 

\begin{introtheorem}\label{thm:intro-cohomology}
Let $n\geq 2$, let $\hess:[n]\to[n]$ be a Hessenberg function, and set
\begin{equation}\label{eq: def D hess}
D(\hess):=\{j\in[n-1]:\hess(j)<\hess(j+1)\}.
\end{equation}
Restriction induces an isomorphism of graded $\mathbb{Q}$-algebras
\[
H^*(\Hess(E_{1,n},\hess);\mathbb{Q})
\cong
\frac{H^*(\Fl(\mathbb{C}^n);\mathbb{Q})}
{\bigl\langle \sigma_{c_j^{(\hess)}}:j\in D(\hess)\bigr\rangle}.
\]
Moreover, every
ideal-generating set of this kernel has at least $|D(\hess)|$ elements, and thus our classes $\{\sigma_{c_j^{(\hess)}}:j\in D(\hess)\}$ form a minimal
generating set. 

\end{introtheorem}

We call the elements of $D(\hess)$ the \emph{corners} of $\hess$; they are the last indices of the maximal intervals on which $\hess$ is constant and less than $n$. It is also worth comparing our description of $J_\hess$ in Theorem~\ref{thm:intro-cohomology} with that of Abe and Crooks. There are $\lvert \mathcal{F}_\hess \rvert = (n-2)!\sum_{j=1}^{n}\bigl(n-\hess(j)\bigr)$ many additive generators for $J_\hess$ in the Abe-Crooks description. Our description of the same ideal has $\lvert D(\hess)\rvert < n$ generators, and these form a minimal set. Thus, a cardinality that is factorial in $n$ has been replaced by one that is linear in $n$.

Our paper is structured as follows.  Section \ref{sec:background}
reviews the necessary
background on flag varieties and Hessenberg varieties.  We also discuss
Schubert polynomials and their connections to divided difference operators.
Our main results are proved in Sections~\ref{sec:mainresults} and~\ref{sec: theorem B}.


\section{Background}\label{sec:background}

We use the flag variety $\Fl(\C^n)$ and the Hessenberg varieties $\Hess(\mx,\hess)$ defined in Section~\ref{sec: intro}, and fix once and for all the matrix
\[
\mx := E_{1,n} = 
\begin{bmatrix}
0 & 0 & \cdots & 0 & 1 \\
0 & 0 & \cdots& 0 & 0 \\
\vdots &  & & & \vdots \\
0 & 0 & \cdots & 0 & 0 
\end{bmatrix} 
\]
with a single $1$ entry in the upper-right corner and $0$'s elsewhere. This matrix $\mx = E_{1,n}$ is in the weight space of the (unique) highest root of $GL_n(\C)$ and belongs to a minimal nilpotent orbit of $GL_n(\C)$. As noted in \cite[Lemma 1]{AbeCrooks}, Hessenberg varieties associated to different matrices in the same orbit are isomorphic, so to study the topology of Hessenberg varieties associated to the minimal nilpotent orbit, it suffices to study $\Hess(\mx=E_{1,n},\hess)$. Thus, in this manuscript we fix $\mx = E_{1,n}$ and allow $\hess$ to vary. 
It is known that, in this case, the Hessenberg variety $\Hess(\mx,\hess)$ is a union of Schubert varieties. 

We now recall the results of Abe and Crooks \cite{AbeCrooks}
regarding the cohomology rings of $\Hess(\mx, \hess)$. These results form the starting point of our investigation. 
A permutation $w \in S_n$ determines a permutation matrix in $GL_n(\C)$ and hence an associated permutation flag, see e.g.~\cite[Section 1.3.3]{Handbook-Hess}. 
Let $T \subset GL(n,\C)$ be the maximal torus of diagonal invertible matrices. Then $T$ acts on $\Fl(\C^n)$, and the $T$-fixed points are precisely the permutation flags corresponding to elements $w \in S_n$~\cite[Exercise 5.63]{Handbook-Hess}. By a slight abuse of notation we will also denote by $w$ the permutation flag, so we write $w \in \Fl(\C^n)$. With this notation we have 
\[
\Fl(\C^n)^T = S_n
\]
and Abe and Crooks prove the following for $\Hess(E_{1,n},\hess)$:

\begin{prop}[{\cite[Proposition 3]{AbeCrooks}}]
\label{prop: fixed points}
Let $\mx = E_{1,n}$ as above and $\hess:[n] \to [n]$ a Hessenberg function. 
Then the set of $T$-fixed points of $\Hess(\mx, \hess)$ is
\begin{equation}\label{eq: hess fixed points}
\Hess(\mx, \hess)^T = \{w\in S_n \, \mid \, w^{-1}(1) \leq \hess(w^{-1}(n)) \}.
\end{equation} 
\end{prop} 

Comparing~\eqref{eq: hess fixed points} with Definition~\ref{def:forbidden}, we see that the $\hess$-forbidden permutations are precisely the permutations that are \emph{not} $T$-fixed points of $\Hess(\mx,\hess)$, that is,
\[
\mathcal{F}_\hess = S_n \setminus \Hess(\mx,\hess)^T.
\]
We emphasize that this description depends on the choice of $\mx = E_{1,n}$.

Next recall that $\Fl(\C^n)$ has an affine paving given by Schubert cells. Let $B$ and $B_-$ denote the Borel subgroups of upper-triangular and lower-triangular matrices in $GL_n(\C)$, respectively. The (opposite) Schubert varieties are defined as $X_-(w) := \overline{B_{-}wB/B} \subset GL_n(\C)/B \cong \Fl(\C^n)$.  Then $X_-(w)$ determines an opposite Schubert class $\sigma_w \in H^{2\ell(w)}(\Fl(\C^n))$. (See \cite[Chapter 1]{Handbook-Hess} for background.) 

The following is a combination of \cite[Corollary 31]{AbeCrooks} and \cite[Theorem 32]{AbeCrooks}. 

\begin{thm}\label{theorem: main AbeCrooks}
Let $\mx = E_{1,n}$ as above and $\hess:[n] \to [n]$ a Hessenberg function. Let $\Hess(\mx,\hess)$ denote the Hessenberg variety associated to $\mx$ and $\hess$. Then 
the linear subspace 
\begin{equation*}
J_\hess := \bigoplus_{w \in \mathcal{F}_\hess} \Q \cdot \sigma_w
\end{equation*} 
is an ideal in $H^*(\Fl(\C^n);\Q)$, and moreover,
\begin{equation*}
H^*(\Hess(\mx, \hess); \Q) \cong H^*(\Fl(\C^n);\Q)/J_\hess.
\end{equation*}
\end{thm}

As mentioned in the introduction, our main technique is to use divided difference operators and Schubert polynomials. We briefly recall the basic definitions and properties we need, following~\cite[Chapters~II and~IV]{Macdonald1991}.  Let $s_i$ denote the simple transposition in $S_n$ which interchanges $i$ and $i+1$ and acts as the identity on other elements.

\begin{defi}\label{def: divided differences}
For $1\leq i<n$, let $s_i$ act on $\Z[x_1,\ldots,x_n]$ by interchanging $x_i$ and $x_{i+1}$. The \emph{divided difference operator} $\partial_i$ on $\Z[x_1,\ldots,x_n]$ is defined by
\begin{equation}\label{eq: definition divided difference i}
\begin{aligned}
    \partial_i f(x_1, x_2, \dots, x_n)
    &:= \frac{f(x_1, x_2, \dots, x_n) - s_i f(x_1, x_2, \dots, x_i,x_{i+1},\dots, x_n)}{x_i - x_{i+1}} \\
    &\phantom{:}= \frac{f(x_1, x_2, \dots, x_n) - f(x_1, \dots, x_{i+1}, x_i, \dots, x_n)}{x_i - x_{i+1}}.
\end{aligned}
\end{equation}
\end{defi}
Note that the divided difference operator satisfies the Leibniz rule
\begin{equation}\label{dd: Leibniz} 
\partial_i(fg)=(\partial_i f)\,g+(s_i f)\,\partial_i g
\end{equation}
for all $f,g\in\Z[x_1,\ldots,x_n]$.

The \emph{Schubert polynomial} $\mathfrak{S}_w(x_1,x_2,\dots, x_n)$ of a permutation $w\in S_n$ is defined recursively using divided difference operations. Specifically, if $w_0 = (n,n-1,\ldots,1)$ is the long word permutation, set 
\begin{equation}\label{dd: w0}
\mathfrak{S}_{w_0}(x_1,x_2,\dots,x_n):=x_1^{n-1}x_2^{n-2}\cdots x_{n-1}, 
\end{equation}\label{dd: length}
and if $s_i\in S_n$ is a simple transposition such that $\ell(ws_i)=\ell(w)+1$, then 
\begin{equation}\label{eq: div diff prop}
\mathfrak{S}_w(x_1,x_2,\dots, x_n)=\partial_i\mathfrak{S}_{ws_i}(x_1,x_2,\dots, x_n).
\end{equation}

The \emph{Lehmer code} of $w\in S_n$ is $L(w)=(L(w)_1,\ldots,L(w)_n)$, where $L(w)_i=|\{k>i : w(k)<w(i)\}|$.  That is,
$L(w)_i$ counts the number of entries to the right of the $i$-th
position in $w$ that are smaller than $w(i)$. We say  $w$ is \emph{dominant} if $L(w)$ is weakly decreasing. We will make use of  the following fact. 

\begin{prop}[{\cite[Equation (4.7)]{Macdonald1991}}]\label{prop: divided difference facts}
Let $w\in S_n$. 
Suppose that $w\in S_n$ is a dominant permutation. Then $\mathfrak{S}_w=x_1^{L(w)_1}\cdots x_n^{L(w)_n}$.
\end{prop}


\section{Proof of Theorem A} \label{sec:mainresults}

In this section we prove Theorem~\ref{thm:intro-polynomial}. The arguments in this section hold over $\Z$, so we work with the polynomial ring $\Z[x_1,\ldots,x_n]$. Here we sketch the outline of our argument. First, we compute the Schubert polynomials $\mathfrak{S}_{c_j^{(\hess)}}$ explicitly (Lemma~\ref{lem: formula Schubert c j hess}). 
Proposition~\ref{prop: theorem intro part 2} is the main technical argument of this section, where we show that the $\mathfrak{S}_{c_j^{(\hess)}}$ generate $\tilde{J}_\hess$. To do so, we define an auxiliary ideal, determined by the positions of $n$ and $1$, and show that it is preserved by the divided difference operators. This is the key step needed for a descending induction argument on Bruhat length. Coupled with the weakly-ascending property of Hessenberg functions, we can express any forbidden Schubert polynomial in terms of those associated to the permutations $c_j^{(\hess)}$, giving the claim. An examination of their leading terms gives the  
Gr\"obner basis statement (Proposition~\ref{prop:grobner}).

The lemma below records facts concerning the
$c_j^{(\hess)}$ which follow directly from Definition~\ref{def:gens}.

\begin{lem}\label{lemma: properties of c j hess}
Let $\hess:[n]\to [n]$ be a Hessenberg function and let $j \in [n]$ satisfy $\hess(j)<n$. Then
\begin{enumerate}
\item $(c_j^{(\hess)})^{-1}(n)=j$,
\item $(c_j^{(\hess)})^{-1}(1)=\hess(j)+1$, and 
\item $c_j^{(\hess)} \in \mathcal{F}_\hess$. 
\end{enumerate}
\end{lem}

\begin{ex}[Running example]\label{ex: running}
    Let $n=6$ and $\hess=(2,3,3,5,5,6)$, and write $c_j:=c_j^{(\hess)}$ throughout this running example. Since $\hess(j)<6$ exactly for $j\leq5$, there are five permutations $c_j$. By Definition~\ref{def:gens}, $c_j$ is the forbidden permutation with $6$ in position $j$, $1$ in position $\hess(j)+1$, and the remaining entries in increasing order. Explicitly, they are 
    $$c_1=(\mathbf{6}\ 2\,\mathbf{1}\,3\,4\,5),\quad
    c_2=(2\,\mathbf{6}\,3\,\mathbf{1}\,4\,5),\quad
    c_3=(2\,3\,\mathbf{6}\,\mathbf{1}\,4\,5),\quad
    c_4=(2\,3\,4\,\mathbf{6}\,5\,\mathbf{1}),\quad
    c_5=(2\,3\,4\,5\,\mathbf{6}\,\mathbf{1}).$$
    \exampleqed
\end{ex}

For $i\in[n]$, we set the notation \[P_i:=x_1x_2\cdots x_i,\] and recall from Theorem~\ref{thm:intro-polynomial} that $H_d(x_1,\ldots,x_j)$ denotes the complete homogeneous symmetric polynomial of degree $d$ in $x_1,\ldots,x_j$. We use three elementary facts. First, $s_jP_i=P_i$ for $j\neq i$, so by Definition~\ref{def: divided differences}, multiplication by $P_i$ commutes with $\partial_j$:
\begin{equation}\label{eq:key3}
\partial_j(P_if)=P_i\,\partial_jf \qquad \textup{ for } \, 1\leq j<n,\ j\neq i,\ f\in\Z[x_1,\ldots,x_n].
\end{equation}
Second, for $d\geq 1$ and $1\leq j<n$, we have~\cite[(3.9)]{Macdonald1991}
\begin{equation}\label{eq: H identity}
\partial_j H_{d}(x_1,\ldots, x_j)=H_{d-1}(x_1,\ldots, x_{j+1}).
\end{equation}
Third, applying~\eqref{eq: H identity} repeatedly to $H_{n-2}(x_1)=x_1^{n-2}$ gives
\begin{equation}\label{eq: complete symmetrics identity}
\partial_{j-1}\cdots\partial_1\,x_1^{n-2}=H_{n-j-1}(x_1,\ldots, x_j) \qquad \textup{ for } \, 1\leq j < n.
\end{equation}

We compute the Schubert polynomials of the permutations $c_j^{(\hess)}$.

\begin{lem}\label{lem: formula Schubert c j hess}
Let $j\in[n]$ with $\hess(j)<n$. Then
\begin{equation*}
\mathfrak{S}_{c_j^{(\hess)}}=P_{\hess(j)}\cdot H_{n-j-1}(x_1,\ldots,x_j).
\end{equation*}
\end{lem}

\begin{proof}
Set $k:=\hess(j)$ and consider the permutation
\[
u=(n,2,\ldots,k,1,k+1,\ldots,n-1).
\]
The number of inversions starting in the first position is $n-1$, the number starting in each of positions $2,\ldots,k$ is $1$, and all remaining inversion counts are zero. These counts form the Lehmer code of $u$ and are weakly decreasing, so $u$ is dominant. By Proposition~\ref{prop: divided difference facts}, we have
\[
\mathfrak{S}_u=x_1^{n-1}x_2\cdots x_k=P_k\cdot x_1^{n-2}.
\]
From the definition of $c_j^{(\hess)}$ we see that $c_j^{(\hess)}=us_1\cdots s_{j-1}$. Each successive right multiplication moves $n$ one position to the right and decreases the length by $1$. Thus repeated application of~\eqref{eq: div diff prop} gives
\[
\mathfrak{S}_{c_j^{(\hess)}}=\partial_{j-1}\cdots\partial_1(P_k\cdot x_1^{n-2}).
\]
Since $j\leq k$, each of $\partial_1,\ldots,\partial_{j-1}$ commutes with multiplication by $P_k$ by~\eqref{eq:key3}, so
\[
\mathfrak{S}_{c_j^{(\hess)}}=P_k\cdot(\partial_{j-1}\cdots\partial_1x_1^{n-2})
=P_k\cdot H_{n-j-1}(x_1,\ldots,x_j),
\]
where the final equality follows from~\eqref{eq: complete symmetrics identity}.
\end{proof}

\begin{ex}[Example \ref{ex: running}, continued]\label{ex: running part 2}
Lemma~\ref{lem: formula Schubert c j hess} gives $\mathfrak{S}_{c_j}=P_{\hess(j)}\cdot H_{n-j-1}(x_1,\ldots,x_j)$ with $\hess(j) <n$.
The table below records the corresponding data for this example. 
\[
\renewcommand{\arraystretch}{1.2}
\begin{array}{c|ccccc}
j & 1 & 2 & 3 & 4 & 5\\ \hline
c_j & 621345 & 263145 & 236145 & 234651 & 234561\\
P_{\hess(j)} & P_2 & P_3 & P_3 & P_5 & P_5\\
H_{6-j-1}(x_1,\ldots,x_j) & x_1^4 & H_3(x_1,x_2) & H_2(x_1,x_2,x_3) & x_1+x_2+x_3+x_4 & 1
\end{array}
\]
The factor $P_{\hess(j)}$ is fixed by the position $\hess(j)+1$ of $1$, and $H_{6-j-1}(x_1,\ldots,x_j)$ by the position $j$ of $6$. Because $\hess(2)=\hess(3)$ and $\hess(4)=\hess(5)$, the pairs $c_2,c_3$ and $c_4,c_5$ share their $P$-factor and differ only in $H_{6-j-1}(x_1,\ldots,x_j)$. The two pairs however behave differently. Since 
$H_{6-5-1}(x_1,\ldots,x_5)=1$, we have
\begin{equation*}
   \mathfrak{S}_{c_4}=P_5(x_1+x_2+x_3+x_4)=(x_1+x_2+x_3+x_4)\mathfrak{S}_{c_5}, 
\end{equation*}
so $c_4$ is redundant already in $\Z[x_1,\ldots,x_6]$. In contrast, $\mathfrak{S}_{c_2}$ is not a multiple of $\mathfrak{S}_{c_3}$. We will revisit this example in Example~\ref{ex: running B}. 
\exampleqed
\end{ex}

Our main result seeks to relate the Schubert polynomials $\mathfrak{S}_w$ for $w\in\mathcal{F}_\hess$ to those of the form $\mathfrak{S}_{c_j^{(\hess)}}$. To do this, we introduce an auxiliary ideal that depends only on the positions of $n$ and $1$.
By symmetry and~\eqref{eq: H identity}, the polynomials 
$H_{n-j-1}(x_1,\ldots,x_j)$ satisfy
\begin{equation}\label{eq: auxiliary f relations}
\begin{array}{rcl}
\partial_i H_{n-j-1}(x_1,\ldots,x_j) &=& 0\quad \textup{ for } i\neq j,\qquad ~~\mbox{and}\\
\partial_j H_{n-j-1}(x_1,\ldots,x_j) & = & H_{n-j-2}(x_1,\ldots,x_{j+1}) \quad \textup{ for } 1 \leq j<n-1.
\end{array}
\end{equation}
For $1\leq p<q\leq n$, define
\begin{equation*}
\mathcal{I}_{p,q}:=
P_{q-1}\langle H_{n-2}(x_1),H_{n-3}(x_1,x_2),\ldots,H_{n-p-1}(x_1,\ldots,x_p) \rangle
\subseteq\Z[x_1,\ldots,x_n].
\end{equation*}

\begin{lem}\label{lemma: schubert poly for forbidden}\label{lem:auxiliary}
Let $w\in S_n$, and suppose that $p=w^{-1}(n)$ and $q=w^{-1}(1)$ satisfy $p<q$. Then
\begin{equation}\label{eq: ideal p q}
\mathfrak{S}_w\in\mathcal{I}_{p,q}.
\end{equation}
\end{lem}
\begin{proof}
We first claim that
\begin{equation}\label{eq: auxiliary stability}
\partial_i\mathcal{I}_{p,q}\subseteq\mathcal{I}_{p,q}
\qquad\textup{if }i\notin\{p,q-1\}.
\end{equation}
Indeed, let $i\notin\{p,q-1\}$, $1\leq j\leq p$, and $g\in\Z[x_1,\ldots,x_n]$. By~\eqref{eq:key3} and~\eqref{dd: Leibniz} we have 
\[
\partial_i(P_{q-1}\cdot g \cdot H_{n-j-1}(x_1,\ldots,x_j))
=P_{q-1}\bigl((\partial_i g)H_{n-j-1}(x_1,\ldots,x_j)
+(s_i g)\partial_i H_{n-j-1}(x_1,\ldots,x_j)\bigr).
\]
Since $i\neq p$, \eqref{eq: auxiliary f relations} shows that $\partial_i H_{n-j-1}(x_1,\ldots,x_j)$ is either $0$ or 
$H_{n-j-2}(x_1,\ldots,x_{j+1})$ with $j+1\leq p$. Hence $\partial_i(P_{q-1}\cdot g \cdot H_{n-j-1}(x_1,\ldots,x_j))\in\mathcal{I}_{p,q}$, from which ~\eqref{eq: auxiliary stability} follows.

We now use descending induction on $\ell(w)$ to prove~\eqref{eq: ideal p q}. The base case is the longest word $w=w_0$. In this case we have $(p,q)=(1,n)$, and
\[
P_{n-1}H_{n-2}(x_1)=x_1^{n-1}x_2\cdots x_{n-1}.
\]
We can then see that the Schubert polynomial $\mathfrak{S}_{w_0}$ is a multiple of $P_{n-1}H_{n-2}(x_1)$ by~\eqref{dd: w0}, so the claim holds in this case. Now suppose $w$ has length $\ell(w)<\ell(w_0)$ and suppose by induction that the claim holds for all permutations with length greater than $k$. Choose an ascent $i$ of $w$, so $w(i)<w(i+1)$, and define $v=ws_i$. Since $w(p)=n$ and $w(q)=1$ by hypothesis, this ascent cannot occur at $p$ or $q-1$. Moreover, if $p'=v^{-1}(n)$ and $q'=v^{-1}(1)$, then $p'\leq p$ and $q'\geq q$. Consequently, $p'<q'$ and 
\[
\mathcal{I}_{p',q'}\subseteq\mathcal{I}_{p,q}.
\]
Since $\ell(v)=\ell(w)+1$, the relation~\eqref{eq: div diff prop} together with the induction hypothesis and~\eqref{eq: auxiliary stability} now give
\[
\mathfrak{S}_w=\partial_i\mathfrak{S}_v
\in\partial_i\mathcal{I}_{p,q}\subseteq\mathcal{I}_{p,q}.
\]
This is what we needed to show.
\end{proof}

We now apply Lemma~\ref{lem:auxiliary} to prove the generation statement. 

\begin{prop}\label{prop: theorem intro part 2}
Let $n\geq 2$, $\hess:[n]\to[n]$ be a Hessenberg function.   
Then, as ideals in $\mathbb{Z}[x_1,\ldots,x_n]$, we have 
\begin{equation}\label{eq: ideal equality polynomial}
\tilde{J}_\hess := \bigl\langle \mathfrak{S}_w : w\in\mathcal{F}_\hess \bigr\rangle
=
\bigl\langle \mathfrak{S}_{c_j^{(\hess)}} : j\in[n],\ \hess(j)<n \bigr\rangle.
\end{equation}
\end{prop}

\begin{proof}
From Lemma~\ref{lemma: properties of c j hess}(3) we know that $c_j^{(\hess)}\in\mathcal{F}_\hess$ for each $j$ with $\hess(j)<n$, so one inclusion is immediate. Thus it suffices to show that if $w\in\mathcal{F}_\hess$, then $\mathfrak{S}_w\in\langle \mathfrak{S}_{c_j^{(\hess)}} : j\in[n],\ \hess(j)<n\rangle$.

Suppose $w\in\mathcal{F}_\hess$, and put $p=w^{-1}(n)$ and $q=w^{-1}(1)$. The forbidden condition and the definition of a Hessenberg function give
\[
p\leq\hess(p)<q.
\]
By Lemma~\ref{lem:auxiliary}, we have $\mathfrak{S}_w \in 
\mathcal{I}_{p,q}$. Therefore, in order to show the claim, it would suffice to show that every generator of $\mathcal{I}_{p,q}$ lies in the RHS of~\eqref{eq: ideal equality polynomial}.  
To see this, let $j \in \Z$ with $1\leq j\leq p$. The weakly ascending condition of $\hess$ implies $\hess(j)\leq\hess(p)<q\leq n$. In particular, $\hess(j)<n$, and Lemma~\ref{lem: formula Schubert c j hess} gives
\begin{equation*}
P_{q-1}H_{n-j-1}(x_1,\ldots,x_j)
=\left(\prod_{r=\hess(j)+1}^{q-1}x_r\right)
P_{\hess(j)}H_{n-j-1}(x_1,\ldots,x_j)
=\left(\prod_{r=\hess(j)+1}^{q-1}x_r\right)
\mathfrak{S}_{c_j^{(\hess)}},
\end{equation*}
where an empty product is understood to be $1$. Thus every generator of $\mathcal{I}_{p,q}$ belongs to the RHS of~\eqref{eq: ideal equality polynomial}, and hence so does $\mathfrak{S}_w$, as required. 
\end{proof}

\begin{rmk}\label{rmk: redundant generators}
The generating set in Theorem~\ref{thm:intro-polynomial} need not be minimal. If $j<n-1$ and $\hess(j)=\hess(n-1)=n-1$, then $\mathfrak{S}_{c_j^{(\hess)}}=(x_1\cdots x_{n-1})\,H_{n-j-1}(x_1,\ldots,x_j)$ is a multiple of $\mathfrak{S}_{c_{n-1}^{(\hess)}}=x_1\cdots x_{n-1}$. 
\end{rmk}

\begin{ex}[Example~\ref{ex: running part 2}, continued]\label{ex: running part 3}
We trace the proof of Theorem~\ref{thm:intro-polynomial}(2) for $w=(362415)$. First note that $6$ is in position $p=2$ and $1$ in position $q=5$, which means $q > \hess(p)=\hess(2)=3$, so $w$ is $\hess$-forbidden. Lemma~\ref{lem:auxiliary} says that  $\mathfrak{S}_w$ lies in $\mathcal{I}_{2,5}=P_4\langle H_4(x_1),H_3(x_1,x_2)\rangle$, and indeed 
\[
\mathfrak{S}_w=x_1^2x_2^2x_3x_4\,(x_1^2+x_1x_2+x_2^2)=P_4\bigl(x_1H_3(x_1,x_2)-H_4(x_1)\bigr).
\]
Since $P_4H_3(x_1,x_2)=x_4\,\mathfrak{S}_{c_2}$ and $P_4H_4(x_1)=x_3x_4\,\mathfrak{S}_{c_1}$, the above formula yields $\mathfrak{S}_w=x_1x_4\,\mathfrak{S}_{c_2}-x_3x_4\,\mathfrak{S}_{c_1}$. A similar argument applies to each of the $|\mathcal{F}_\hess|=4!\cdot 12=288$ forbidden permutations. Thus, in this example, we have replaced the $288$ generators with $5$ generators. 
\exampleqed
\end{ex}

Our Theorem \ref{thm:intro-polynomial} also makes a claim about Gr\"obner bases. Our next result addresses this.

\begin{prop}\label{prop:grobner}
Let $n\geq 2$ and $\hess:[n]\to[n]$ be a Hessenberg function. Consider the lexicographic order $x_n>\cdots>x_1$ on $\Z[x_1,\ldots,x_n]$. Then the polynomials $\mathfrak{S}_{c_j^{(\hess)}}$, for $j\in[n]$ with $\hess(j)<n$, form a monic Gr\"obner basis of $\tilde{J}_\hess$, with leading monomials
\[
L_j := P_{\hess(j)} \cdot x_j^{n-j-1}.
\]
\end{prop}

\begin{proof}
For the given lexicographic order, it is immediate from the definition of $H_{n-j-1}(x_1,\ldots,x_j)$ that its leading term is $x_j^{n-j-1}$, which has coefficient $1$. Multiplying by the monomial $P_{\hess(j)}$ then yields that the leading monomial of $P_{\hess(j)}\cdot H_{n-j-1}(x_1,\cdots,x_j)$ is $L_j = P_{\hess(j)} \cdot x_j^{n-j-1}$, which also has coefficient $1$. 
In particular, each $\mathfrak{S}_{c_j^{(\hess)}}$ is monic, and it remains to check the $S$-polynomials.

Let $j_1,j_2 \in [n]$ with $\hess(j_1)<n, \hess(j_2)<n$ and assume $j_1<j_2$. For simplicity of notation, set $g_1 := \mathfrak{S}_{c_{j_1}^{(\hess)}}$ and $g_2 := \mathfrak{S}_{c_{j_2}^{(\hess)}}$. We now compute the $S$-polynomial $S(g_1,g_2)$. 
From the assumption $j_1< j_2$ we have $\hess(j_1) \leq \hess(j_2)$, which means 
\[
\mathbf{m}:=\gcd(L_{j_1},L_{j_2})=P_{\hess(j_1)}.
\]
From the formula for $\mathfrak{S}_{c_j^{(\hess)}}$ we also know that the monomial $\mathbf{m}$ divides both $g_1$ and $g_2$.
In particular, if $T_{j_1},T_{j_2}$ denote the remaining terms of $g_1,g_2$ so that 
\[
g_1=L_{j_1}+T_{j_1},\qquad g_2=L_{j_2}+T_{j_2},
\]
then we also have $\mathbf{m}\mid T_{j_1}$ and $\mathbf{m}\mid T_{j_2}$.
Let
\[
M:=\operatorname{lcm}(L_{j_1},L_{j_2})=\frac{L_{j_1}L_{j_2}}{\mathbf{m}}
\]
denote the lcm of $L_{j_1},L_{j_2}$. 
Since $g_1$ and $g_2$ are monic, the $S$-polynomial is
\begin{align*}
S(g_1,g_2)
&=\frac{M}{L_{j_1}}g_1-\frac{M}{L_{j_2}}g_2\\
&=\frac{L_{j_2}}{\mathbf{m}}(L_{j_1}+T_{j_1})
  -\frac{L_{j_1}}{\mathbf{m}}(L_{j_2}+T_{j_2})\\
&=\frac{L_{j_2}T_{j_1}-L_{j_1}T_{j_2}}{\mathbf{m}}\\
&=\frac{T_{j_1}(g_2-T_{j_2})-T_{j_2}(g_1-T_{j_1})}{\mathbf{m}} =\frac{T_{j_1}}{\mathbf{m}}g_{2}-\frac{T_{j_2}}{\mathbf{m}}g_1.
\end{align*}
Both coefficients in the last expression belong to
$\mathbb{Z}[x_1,\ldots,x_n]$.

Every monomial in $T_{j_1}$ is strictly smaller in the given lexicographic order than $L_{j_1}$,
and likewise for $T_{j_2}$. Thus, omitting any zero summands,
\begin{align*}
\operatorname{LM}\!\left(\frac{T_{j_1}}{\mathbf{m}}g_2\right)
&=\frac{\operatorname{LM}(T_{j_1})L_{j_2}}{\mathbf{m}}
 <\frac{L_{j_1}L_{j_2}}{\mathbf{m}}=M,\\
\operatorname{LM}\!\left(\frac{T_{j_2}}{\mathbf{m}}g_1\right)
&=\frac{\operatorname{LM}(T_{j_2})L_{j_1}}{\mathbf{m}}
 <\frac{L_{j_1}L_{j_2}}{\mathbf{m}}=M.
\end{align*}
Hence every $S$-polynomial has a representation in the
generators in which each summand has leading monomial
strictly smaller than the corresponding least common
multiple. This is an LCM representation in the sense of \cite[Chapter~2, \S~9, Definition~5]{CoxLittleOShea2015}. A refinement of Buchberger's criterion in terms of LCM representations \cite[Chapter~2, \S~9, Theorem~6]{CoxLittleOShea2015} implies 
that the $g_j$ form a Gr\"obner basis over $\mathbb{Q}$. To obtain the integral statement, take
$f\in\tilde{J}_\hess\subseteq\Z[x_1,\ldots,x_n]$.
Division by the monic polynomials $g_j$ uses only
integer coefficients. Its remainder is also a
remainder for the same division over $\Q$, where it
is zero by the Gr\"obner-basis property.
Thus the $g_j$ form a monic Gr\"obner basis over $\Z$.
\end{proof}

\begin{rmk}
The exponent vectors of the leading monomials in Proposition~\ref{prop:grobner}  are the Lehmer codes of the permutations $c_j^{(\hess)}$, as defined before Proposition~\ref{prop: divided difference facts}. Indeed, 
\[
L(c_j^{(\hess)}) = (\underbrace{1,\ldots,1}_{j-1},\, n-j,\, \underbrace{1,\ldots,1}_{\hess(j)-j},\, \underbrace{0,\ldots,0}_{n-\hess(j)}),
\]
which is precisely the exponent vector of $L_j=P_{\hess(j)}\cdot x_j^{n-j-1}$.
These leading terms can also be interpreted as the monomials associated to the bottom pipe dreams of the tuple $c_j^{(\hess)}$, where crosses in row $i$ occupy the first $L(c_j^{(\hess)})_i$ columns; see~\cite{BergeronBilley1993}.
\end{rmk}

\begin{rmk}\label{remark: redundancy 2}
We noted in Remark~\ref{rmk: redundant generators} that if $\hess(n-1)=n-1$, then the generators corresponding to $j<n-1$ with $\hess(j)=\hess(n-1)=n-1$ are multiples of $\mathfrak{S}_{c_{n-1}^{(\hess)}}$ and hence are not needed. 
These are the only divisibilities among the leading monomials of the $\mathfrak{S}_{c_j^{(\hess)}}$ appearing in Proposition~\ref{prop:grobner}, so by discarding the generators corresponding to these $j$  we still obtain a minimal Gr\"obner basis. 
\end{rmk}

\begin{theoremproof}{thm:intro-polynomial}
Part (1) of the theorem is in Lemma~\ref{lem: formula Schubert c j hess}, part (2) is Proposition~\ref{prop: theorem intro part 2}, and Proposition~\ref{prop:grobner} establishes the Gr\"obner basis claim in (3). This completes the proof.
\end{theoremproof}

\section{Proof of Theorem B}\label{sec: theorem B}

In this section, we prove Theorem~\ref{thm:intro-cohomology}. We first sketch our argument. The ring $H^*(\Fl(\C^n);\Q)$ is a quotient of the polynomial ring $\Q[x_1,\cdots,x_n]$ by the ideal $\langle \mathbf{e}_1,\cdots,\mathbf{e}_n\rangle$.
Our key observation is that if $j \in [n-1]$ is such that $\hess(j)=\hess(j+1)<n$, then the Schubert class $\sigma_{c_{j}^{(\hess})}$ is a multiple of $\sigma_{c_{j+1}^{(\hess)}}$ (cf. Proposition~\ref{prop:jumps are enough}). This allows us to reduce to the corner generators. Next, we
identify the $c_j^{(\hess)}$ for $j \in D(\hess)$ as the Bruhat-minimal
elements of the forbidden set $\mathcal{F}_\hess$. The Bruhat support of Schubert
multiplication shows that these classes remain linearly independent
modulo positive-degree multiples of the kernel, establishing the minimality and cardinality claims. 
We use $\Q$ coefficients, following \cite{AbeCrooks}.

For a partially ordered set $(\mathcal{S},\leq)$, we say that a subset $\mathcal{T}$ of $\mathcal{S}$ is a \emph{down-set} if, whenever $w \in \mathcal{T}$ and $u \leq w$, then $u \in \mathcal{T}$. Similarly we say that $\mathcal{U}$ is an \emph{up-set} if, whenever $w \in \mathcal{U}$ and $u \geq w$, then $u \in \mathcal{U}$. Note that the complement of a down-set is an up-set, and vice versa.

\begin{lem}\label{lem:filter}
Let $n \geq 2$ and $\hess: [n] \to [n]$ be a Hessenberg function. Let $\mx = E_{1,n}$ and let $\Hess(\mx,\hess)$ be the Hessenberg variety corresponding to $\mx$ and $\hess$. Then $\Hess(\mx,\hess)^T$ is a down-set with respect to the Bruhat order, and $\mathcal{F}_\hess = S_n \setminus \Hess(\mx,\hess)^T$ is an up-set with respect to the Bruhat order.  
\end{lem}

We will use the \emph{tableau criterion} \cite[Theorem~2.6.3]{BjornerBrenti2005} for the Bruhat order 
on $S_n$. This is defined by creating a tableau for a permutation
$u = (u_1,u_2,\ldots,u_n)$ such that in row $1 \leq i \leq n-1$,
we place the first $n-i$ elements of $u$ in ascending order.  
If $u,w \in S_n$, 
then $u \leq w$ with respect to the Bruhat order if and only if every entry
of the tableau of $u$ is less than or equal to the corresponding entry in
the tableau of $w$. For example, if $u = (21453)$ and $w=(32541)$, then
the tableaux are 
\[
\begin{array}{cccc}
1&2&4&5 \\
1&2&4& \\
1&2& & \\
2 &&&
\end{array}~~\mbox{and}\hspace{.5cm}
\begin{array}{cccc}
2&3&4&5 \\
2&3&5& \\
2&3& & \\
3 &&&
\end{array}
\]
Then $u \leq w$ in Bruhat order since each entry in the tableau of $u$ is less than or equal to the corresponding entry in the tableau of $w$.

\begin{proof}
Suppose $u, w \in S_n$.  Suppose also that $u\leq w$ and $w \in \Hess(\mx,\hess)^T$. From the tableau criterion for Bruhat order, if $u \leq w$, then $u^{-1}(1)\le w^{-1}(1)$ and $u^{-1}(n)\geq w^{-1}(n)$.  Since $w \in \Hess(\mx,\hess)^T$ and $\hess$ is weakly increasing, we conclude that 
\[u^{-1}(1)\leq w^{-1}(1)\leq \hess(w^{-1}(n))\leq\hess(u^{-1}(n))\]
which means that $u\in\Hess(\mx,\hess)^T$, as desired. The second claim follows immediately. 
\end{proof}

Recall from Section~\ref{sec: intro} that the set $D(\hess)$ of corners is $\{ j\in[n-1]: \hess(j)<\hess(j+1)\}$. We begin with a characterization of when a permutation lies above $c_j^{(\hess)}$ in Bruhat order, in terms of its one-line notation.

\begin{lem}\label{lemma: Bruhat bigger than c j hess}
In the setting of Lemma~\ref{lem:filter}, let $j \in [n-1]$ and suppose $\hess(j)<n$, and let $c_j^{(\hess)}$ denote the permutation of Definition~\ref{def:gens}. Let $w\in S_n$. Then:
\[
c_j^{(\hess)} \leq w \Leftrightarrow w^{-1}(n) \leq j \, \textup{ and } \, w^{-1}(1) \geq \hess(j) +1. 
\]
\end{lem}

\begin{proof} 
First suppose that $c_j^{(\hess)} \leq w$. By Lemma~\ref{lemma: properties of c j hess} we know that $n$ appears in the $j$-th position and $1$ appears in the $(\hess(j)+1)$-th position of the one-line notation of $c_j^{(\hess)}$. By the tableau criterion for Bruhat order, it follows that $n$ must appear in one of the first $j$ positions of $w$, i.e., $w^{-1}(n) \leq j$, and similarly that $1$ must appear in the $(\hess(j)+1)$-th position or further to its right, i.e., $w^{-1}(1) \geq \hess(j)+1$.

For the reverse direction, suppose now that $w^{-1}(n) \leq j$ and $w^{-1}(1)\geq \hess(j)+1$. By the tableau criterion, it suffices to show that for each $1\leq a\leq n-1$, the $i$-th smallest of the first $a$ entries of $c_j^{(\hess)}$ is at most the $i$-th smallest of the first $a$ entries of $w$, for every $i\leq a$. By Definition~\ref{def:gens}, the first $a$ entries of $c_j^{(\hess)}$ form the set
\[
\{2,\ldots,a+1\}\ \ \textup{ for } \, a<j,\quad
\{2,\ldots,a\}\cup\{n\}\ \ \textup{ for} \, j\leq a\leq\hess(j),\quad
\{1,\ldots,a-1\}\cup\{n\}\ \ \textup{ for } a>\hess(j).
\]
If $a\leq\hess(j)$, then $1$ is not among the first $a$ entries of $w$, so the $i$-th smallest of them is at least $i+1$. If $a\geq j$, then $n$ is among the first $a$ entries of $w$, so the largest of them is $n$. Together with the fact that the $i$-th smallest of any $a$ distinct positive integers is at least $i$, these two observations give the required inequalities in each of the three cases, so $c_j^{(\hess)} \leq w$.
\end{proof}

\begin{lem}\label{lem:minimal-jumps}
In the setting of Lemma~\ref{lem:filter}, the Bruhat-minimal elements of $\mathcal{F}_\hess$ are precisely the $c_j^{(\hess)}$ where $j \in D(\hess)$ is a corner of $\hess$.
\end{lem}

\begin{proof}
Suppose $w\in\mathcal{F}_\hess$. Let $j:=w^{-1}(n)$ be the position of $n$ in the one-line notation of $w$. Since $w$ is in $\mathcal{F}_\hess$, this means $w^{-1}(1) > \hess(j)$. Lemma~\ref{lemma: Bruhat bigger than c j hess} then implies $c_j^{(\hess)} \leq w$.
Thus every element of $\mathcal{F}_\hess$ lies above some member of the family
$\{c_j^{(\hess)}: j \in [n-1], \hess(j)<n\}$, and this family is contained in $\mathcal{F}_\hess$ by Lemma~\ref{lemma: properties of c j hess}(3). Hence the Bruhat-minimal elements of $\mathcal{F}_\hess$ are exactly the minimal elements of this family.
Next suppose that $j,j' \in [n-1]$ and $\hess(j)<n$ and $\hess(j')<n$. By Lemma~\ref{lemma: Bruhat bigger than c j hess} we have  
\[
c_j^{(\hess)}\le c_{j'}^{(\hess)}
\quad\Longleftrightarrow\quad
j'\le j\ \text{ and }\ \hess(j') \geq \hess(j).
\]
By the definition of Hessenberg functions, if $j'\leq j$, then $\hess(j') \leq \hess(j)$ so in fact we have 
\begin{equation}\label{eq:bruhat order on cj}
  c_j^{(\hess)}\le c_{j'}^{(\hess)}
\quad\Longleftrightarrow\quad
j'\le j\ \text{ and }\ \hess(j') =\hess(j).  
\end{equation}

Thus, two indices $j$ and $j'$ with distinct values of $\hess$ have corresponding $c_j^{(\hess)}$ and $c_{j'}^{(\hess)}$ which are
incomparable in Bruhat order, and for a collection of (consecutive) indices with the same value of $\hess$, it is the largest value of the index which yields the Bruhat-smallest $c_j^{(\hess)}$. 
The largest such indices are precisely the corners $j\in D(\hess)$. This proves the lemma.
\end{proof}

\begin{ex}[Example~\ref{ex: running}, continued] 
By Equation~\ref{eq:bruhat order on cj},  $c_j\leq c_{j'}$ if and only if $j'\leq j$ and $\hess(j')=\hess(j)$.  The only such pairs with $j\neq j'$ are $c_3<c_2$ and $c_5<c_4$. For the first pair, the tableaux
\[
c_3=(236145):
\begin{array}{ccccc}
1&2&3&4&6\\
1&2&3&6\\ 
2&3&6\\ 
2&\mathbf{3}\\ 
2
\end{array}
\qquad\qquad
c_2=(263145):
\begin{array}{ccccc}
1&2&3&4&6\\ 
1&2&3&6\\ 
2&3&6\\ 
2&\mathbf{6}\\
2
\end{array}
\]
differ only in the fourth row. For example, $c_1\not\leq c_3$ because $c_3^{-1}(6)=3>1$, and $c_3\not\leq c_1$ because $c_1^{-1}(1)=3<\hess(3)+1$. All other pairs are incomparable. Hence the Bruhat-minimal elements of $\mathcal F_\hess$ are $c_1,c_3,c_5$, indexed by the corners $D(\hess)=\{1,3,5\}$, in accordance with  Lemma~\ref{lem:minimal-jumps}.
\exampleqed
\end{ex}

We saw in the above lemma that if $\hess(j)=\hess(j+1)<n$ then $c_j^{(\hess)} \geq c_{j+1}^{(\hess)}$. There is a corresponding statement about the associated Schubert polynomials. Recall that $\mathbf{e}_1,\ldots,\mathbf{e}_n$ denote the elementary symmetric polynomials in the variables $x_1,\ldots,x_n$. Before proceeding, it is useful to recall the following relation between complete homogeneous symmetric polynomials and the elementary symmetric polynomials.  This is standard but we include it here for the convenience of the reader. Let $\mathbf{H}_A(t)$ denote the generating function for the complete homogeneous symmetric functions $H_d$ on a set of variables $A$, and similarly $\mathbf{E}_A(t)$ for the elementary symmetric functions $\mathbf{e}_d$. 
(So for example if $A = \{x_1,\ldots,x_\ell\}$ for some $\ell>0$,
then $\mathbf{H}_A(t) = \prod_{i=1}^\ell \frac{1}{1-x_i t}$ and $\mathbf{E}_A(t) = \prod_{i=1}^{\ell}(1+x_it)$.) 

\begin{lem}\label{lemma: lies in ideal}
Let $n \geq 2$. Let $j, \ell$ be integers such that $1 \leq j \leq n-2$ and $j<\ell<n$. Then $x_1 \cdots x_\ell \cdot H_{n-j-1}(x_1,\ldots,x_{j+1})$ lies in the ideal $\langle \mathbf{e}_1,\ldots,\mathbf{e}_n\rangle$.
\end{lem}

\begin{proof} 
The argument uses well-known identities in symmetric functions; we sketch the argument. By the definitions above, we have the identity on generating functions $\mathbf{H}_A(t)\mathbf{E}_A(-t)=1$. From this we immediately obtain $\mathbf{H}_{\{x_1,\ldots,x_{j+1}\}}(t)\mathbf{E}_{\{x_1,\ldots,x_n\}}(-t) = \mathbf{E}_{\{x_{j+2},\ldots,x_n\}}(-t)$. Now compare the coefficient of $t^{n-j-1}$ on both sides to obtain 
\[
\sum_{k=0}^{n-j-1}(-1)^k \mathbf{e}_k \cdot
H_{n-j-1-k}(x_1,\ldots,x_{j+1})
=(-1)^{n-j-1} x_{j+2}\cdots x_n.
\]
Rearranging to isolate the term corresponding to $k=0$ on the LHS, multiplying by $x_1\cdots x_\ell$, and using $x_1\cdots x_\ell \cdot x_{j+2}\cdots x_n=x_{j+2}\cdots x_\ell \cdot \mathbf{e}_n$ (where we use $\ell>j$, and we interpret the notation $x_{j+2}\cdots x_\ell$ as equal to $1$ if $j+2>\ell$), gives
\begin{equation*}
    \begin{split}
x_1 \cdots x_\ell \cdot H_{n-j-1}(x_1,\ldots,x_{j+1}) &= (-1)^{n-j-1} x_{j+2}\cdots x_\ell \cdot \mathbf{e}_n \\
& + \sum_{k=0}^{n-j-2} (-1)^k \cdot x_1\cdots x_\ell \cdot \mathbf{e}_{k+1} \cdot H_{n-j-2-k}(x_1,\ldots,x_{j+1}).
\end{split}
\end{equation*}
Every term on the right-hand side lies in 
the ideal $\langle \mathbf{e}_1,\ldots,\mathbf{e}_n\rangle$, so the conclusion
holds.
\end{proof}

We can now prove the last remaining step for the proof of Theorem~\ref{thm:intro-cohomology}, which is related to non-corners
of $\hess$.

\begin{prop}\label{prop:jumps are enough}
In the setting of Lemma~\ref{lem:filter}, suppose $j \in [n-1]$
is such that $\hess(j)=\hess(j+1)<n$. Then
\[
\mathfrak{S}_{c_j^{(\hess)}}
\equiv
-x_{j+1}\mathfrak{S}_{c_{j+1}^{(\hess)}}
\pmod{\langle \mathbf{e}_1,\ldots,\mathbf{e}_n\rangle}.
\]
\end{prop}

\begin{proof}
We begin by noting that the assumption $\hess(j+1)<n$ implies that $j+1<n$ and hence $j \leq n-2$, so $n-j-2 \geq 0$. Moreover, since $\hess(j+1)\geq j+1$ by assumption on Hessenberg functions, we know $j < \hess(j)=\hess(j+1)$. 
Next observe that, for any positive degree $d$, 
we have the relation $H_d(x_1,\dots,x_j)=H_d(x_1,\dots,x_{j+1})-x_{j+1}H_{d-1}(x_1,\dots,x_{j+1})$. Using this with $d=n-j-1\geq 1$,  we have
\begin{equation*}
    \begin{split}
    \mathfrak{S}_{c_j^{(\hess)}}
&=P_{\hess(j)}H_{n-j-1}(x_1,\ldots,x_j) \quad \textup{ by Lemma~\ref{lem: formula Schubert c j hess}} \\
&=P_{\hess(j)}H_{n-j-1}(x_1,\dots,x_{j+1})-x_{j+1}P_{\hess(j)}H_{n-j-2}(x_1,\ldots,x_{j+1})\\
&=P_{\hess(j)}H_{n-j-1}(x_1,\dots,x_{j+1})-x_{j+1}P_{\hess(j+1)}
H_{n-j-2}(x_1,\ldots,x_{j+1})
\, \textup{ since $\hess(j)=\hess(j+1)$ } \\
&=P_{\hess(j)}H_{n-j-1}(x_1,\dots,x_{j+1})-x_{j+1}\,\mathfrak{S}_{c_{j+1}^{(\hess)}} \, \textup{ by Lemma~\ref{lem: formula Schubert c j hess}.}
\end{split} 
\end{equation*}
To conclude the argument it suffices to see that $P_{\hess(j)} \cdot H_{n-j-1}(x_1,\ldots,x_{j+1})$ is contained in the ideal $\langle \mathbf{e}_1,\ldots,\mathbf{e}_n\rangle$. Since $\hess(j)>j$ as noted above, this follows from Lemma~\ref{lemma: lies in ideal}.  This concludes the proof. 
\end{proof}

\begin{ex}[Example~\ref{ex: running part 2}, continued]\label{ex: running B}
For the running example, Proposition~\ref{prop:jumps are enough} applies to $j=2$ and $j=4$, since $\hess=(2,3,3,5,5,6)$ and thus $\hess(2)=\hess(3)=3<6$ and $\hess(4)=\hess(5)=5<6$. For $j=2$, we follow the argument in the proof and observe 
$H_3(x_1,x_2)=H_3(x_1,x_2,x_3)-x_3H_2(x_1,x_2,x_3)$, which gives $\mathfrak{S}_{c_2}=P_3H_3(x_1,x_2,x_3)-x_3\mathfrak{S}_{c_3}$. The first term lies in $\langle\mathbf{e}_1,\ldots,\mathbf{e}_6\rangle$ by Lemma~\ref{lemma: lies in ideal} applied to the choice $(j,\ell)=(2,3)$. Indeed, writing $H_d=H_d(x_1,x_2,x_3)$, we get
\[
P_3H_3=P_3\bigl(\mathbf{e}_1H_2-\mathbf{e}_2H_1+\mathbf{e}_3\bigr)-\mathbf{e}_6,
\]
which visibly lies in the ideal $\langle \mathbf{e}_1,\cdots,\mathbf{e}_6\rangle$. 
Hence $\mathfrak{S}_{c_2}\equiv -x_3\,\mathfrak{S}_{c_3}\pmod{\langle\mathbf{e}_1,\ldots,\mathbf{e}_6\rangle}$.

For $j=4$ we have $H_1(x_1,x_2,x_3,x_4)=(x_1+\cdots+x_5)-x_5 = (x_1+\cdots+x_6)-x_6 - x_5$ and we also know $\mathfrak{S}_{c_5}=P_5$ and $\mathfrak{S}_{c_4}=P_5 \cdot H_1(x_1,x_2,x_3,x_4)$. Observe also that $P_5 \cdot x_6 = P_6 = \mathbf{e}_6$. Putting these together yields
\[
\mathfrak{S}_{c_4}
=P_5\,(\mathbf{e}_1 - x_6)-\,x_5\,P_5
 = P_5 \mathbf{e}_1 - \mathbf{e}_6 - x_5 \mathfrak{S}_{c_5} \equiv -x_5\,\mathfrak{S}_{c_5}\pmod{\langle\mathbf{e}_1,\ldots,\mathbf{e}_6\rangle}.
\]
The above computation follows the argument made in the proof of Proposition~\ref{prop:jumps are enough}. However, as already observed in Remark~\ref{rmk: redundant generators}, 
in the case of $j=4$ we can view the relation in a different way. Since $\mathfrak{S}_{c_5}=P_5$ and $\mathfrak{S}_{c_4}=P_5 \cdot H_1(x_1,x_2,x_3,x_4)$, we can simply observe that $\mathfrak{S}_{c_4}$ is a multiple of $\mathfrak{S}_{c_5}$, and hence is a redundant generator.

The remaining generators are indexed by the corners $D(\hess)=\{1,3,5\}$, and $c_1,c_3,c_5$ are exactly the Bruhat-minimal forbidden permutations (Lemma~\ref{lem:minimal-jumps}). Theorem~\ref{thm:intro-cohomology} therefore gives
\[
H^*(\Hess(E_{1,6},\hess);\Q)\cong H^*(\Fl(\C^6);\Q)\big/\langle\sigma_{c_1},\sigma_{c_3},\sigma_{c_5}\rangle,
\]
and no generating set of this kernel has fewer than three elements, compared with the $288$ classes of the Abe--Crooks additive basis.
\exampleqed
\end{ex}

We now show that the classes indexed by the Bruhat-minimal forbidden permutations generate $J_\hess$, thus proving the second 
main result of this paper.

\begin{theoremproof}{thm:intro-cohomology}
We begin by recalling what we know in the context of the polynomial ring $\Q[x_1,\ldots,x_n]$, building on the fact that $H^*(\Fl(\C^n);\Q) \cong \Q[x_1,\ldots,x_n]/\langle \mathbf{e}_1,\ldots,\mathbf{e}_n\rangle$. First note that the proof of Theorem~\ref{thm:intro-polynomial} also holds verbatim with $\Q$ coefficients in place of $\Z$ coefficients. With this said, from here onwards we interpret all statements to be in $\Q[x_1,\cdots,x_n]$. By Theorem~\ref{thm:intro-polynomial}(2), the ideal $\tilde{J}_\hess = \langle \mathfrak{S}_w: w \in \mathcal{F}_\hess\rangle$ is
generated by the polynomials $\mathfrak{S}_{c_j^{(\hess)}}$ with
$\hess(j)<n$. Since $\sigma_w$ is represented by $\mathfrak{S}_w$ and, by Theorem~\ref{theorem: main AbeCrooks}, $J_\hess$ is the ideal generated by $\{\sigma_w : w\in\mathcal{F}_\hess\}$, it follows that $J_\hess$ is generated by the classes $\sigma_{c_j^{(\hess)}}$ with $\hess(j)<n$. Moreover, Proposition~\ref{prop:jumps are enough} gives
\[
\sigma_{c_j^{(\hess)}}
= -x_{j+1}\sigma_{c_{j+1}^{(\hess)}} \in H^{*}(\Fl(\C^n);\Q)
\qquad\text{if }\hess(j)=\hess(j+1)<n
\]
where by slight abuse of notation, $x_{j+1}$ also denotes its image in $H^*(\Fl(\C^n);\Q)$. By the definition of the set $D(\hess)$ given in~\eqref{eq: def D hess} it follows that 
\[
J_\hess=
\bigl\langle \sigma_{c_j^{(\hess)}}:j\in D(\hess)\bigr\rangle.
\]

We next prove that every ideal-generating set has at
least $|D(\hess)|$ elements. First recall that the Schubert classes form a $\Q$-additive basis for $H^*(\Fl(\C^n);\Q)$ and that $\{\sigma_w: w \in \mathcal{F}_\hess\}$ form a $\Q$-additive basis for $J_\hess$. Let $\mathfrak{m}$ denote the ideal of positive-degree elements in $H^*(\Fl(\C^n);\Q)$. The subspace $\mathfrak{m} \cdot J_\hess$ is spanned by products of the form $\sigma_u \sigma_v$ for $\ell(u)>0$ and $v \in \mathcal{F}_\hess$. By iterating the variable-multiplication formula~\cite[(4.15$'$)]{Macdonald1991}, and using commutativity, we know $\sigma_u \sigma_v = \sum_{w \geq u,v} c^w_{u,v} \sigma_w$ for $c^w_{u,v} \in \Z$ and if $c^w_{u,v} \neq 0$, then $\ell(w)=\ell(u)+\ell(v)$. In particular, any $\sigma_w$ appearing on the RHS must have $w \geq v$, and since $\ell(u)>0$ this implies $\ell(w)>\ell(v)$ and hence the inequality is strict, i.e., $w>v$. Since $\{\sigma_z: z \in \mathcal{F}_\hess\}$ is an additive basis of $J_\hess$ by Theorem~\ref{theorem: main AbeCrooks} we know $w \in \mathcal{F}_\hess$, and by Lemma~\ref{lem:minimal-jumps} and the fact that $w>v$ we conclude that $w$ cannot be a minimal Bruhat element. It follows that, for any element $f$ of $\mathfrak{m}J_\hess$, the basis expansion 
of $f$ in terms of the basis $\{\sigma_z: z \in \mathcal{F}_\hess\}$ has a $0$ coefficient in front of all the Bruhat-minimal elements $\{\sigma_{c_j^{(\hess)}}: j \in D(\hess)\}$. This in turn implies that the image of $\{\sigma_{c_j^{(\hess)}}: j \in D(\hess)\}$ is linearly independent in the $\Q$-vector space $J_\hess/\mathfrak{m}J_\hess$. On the other hand, they must also span this vector space, since we just saw above that they generate $J_\hess$. It follows that 
\[
\dim_\Q(J_\hess/\mathfrak mJ_\hess)=|D(\hess)|.
\]
Observe that for any set $S$ that generates the ideal $J_\hess$, the image of $S$ in $J_\hess/\mathfrak{m}J_\hess$ must $\Q$-span this quotient, since $H^*(\Fl(\C^n);\Q)/\mathfrak{m} \cong \Q$. Hence every such $S$ must have at least $\lvert D(\hess)\rvert$ many elements. The set $\{\sigma_{c_j^{(\hess)}}: j \in D(\hess)\}$ has precisely $\lvert D(\hess)\rvert$ elements and hence it is minimal, i.e., there is no proper subset which generates $J_\hess$.

Finally, the asserted presentation of the cohomology
ring follows from Theorem~\ref{theorem: main AbeCrooks}.
\end{theoremproof}

\subsection*{Acknowledgments}
We thank Peter Crooks for helpful conversations, interest in this work, and ongoing encouragement. We additionally thank Colleen Robichaux, Anna Weigandt, and Oliver Pechenik for sharing their expertise on Schubert polynomials and pipe dreams in the early stages of our investigations. The set of generators $\{c_j^{(\hess)}\}$ as discussed in this manuscript was first found experimentally through computations with SAGE, as part of Nardone's NSERC USRA project in Summer 2026 while working under the 
supervision of Harada.  
GG was funded by SNF Postdoc.Mobility grant P500PT\textunderscore225436.
MH was supported by a Canada Research Chair Award (Tier 1) and NSERC Discovery Grant 2019-06567. 
JR was supported by NSERC Discovery Grant 2023-04800. 
AVT was supported by NSERC Discovery Grant 2024--05299. 

\subsection*{Declaration of AI Use}  
The authors made use of Claude Opus 5.5 and GPT Astra 6 as
part of their investigation of the minimal generators.  
The authors reviewed  all AI output and  wrote
the paper; they take full 
responsibility for the content of this paper.

\bibliographystyle{amsplain}
\bibliography{ref}

\end{document}